\documentclass[12pt]{article}

\usepackage{amssymb}
\usepackage{amsfonts}
\usepackage{latexsym}
\usepackage{graphicx}
\usepackage{amsmath}
\usepackage{amsthm}
\usepackage{lscape}

\newcommand{\aut}{\mathrm{aut}\,}

\newcommand{\beqn}{\begin{eqnarray*}}
\newcommand{\eeqn}{\end{eqnarray*}}

\newcommand{\hd}{\mathcal{H}}
\newcommand{\tl}{\mathcal{T}}

\newcommand{\Z}{\mathbb{Z}}
\newcommand{\R}{\mathbb{R}}
\newcommand{\C}{\mathbb{C}}
\newcommand{\D}{\mathbb{D}}

\newcommand{\ES}{\mathbb{S}}

\newcommand{\ONE}{\mathbf{1}}  

\renewcommand{\AA}{{\mathcal A}} 

\newcommand{\CC}{{\mathcal C}}

\newcommand{\GG}{{\mathcal G}}

\newcommand{\KK}{{\mathcal K}}
\newcommand{\WW}{{\mathcal W}}

\newcommand{\Matrix}[1]{\ensuremath{\left[\begin{array}{ccccccccccccccccccccccccr} #1 \end{array}\right]}}

\newtheorem{theorem}{Theorem}[section]
\newtheorem{lemma}[theorem]{Lemma}
\newtheorem{definition}[theorem]{Definition}

\newtheorem{proposition}[theorem]{Proposition}
\newtheorem{remark}[theorem]{Remark}
\newtheorem{example}[theorem]{Example}

\newcommand{\Fix}{\mathrm{Fix}}

\title{Balanced Colorings and Bifurcations\\in Rivalry and Opinion Networks}
\author{Ian Stewart \\Mathematics Institute\\University of Warwick\\Coventry CV4 7AL \\ United Kingdom} 
\date{\today}

\begin{document}
\maketitle

\begin{abstract}
Balanced colorings of networks classify robust synchrony patterns --- those that are 
defined by subspaces that are flow-invariant for all admissible ODEs. In
symmetric networks the obvious balanced colorings are orbit colorings,
where colors correspond to orbits of a subgroup of the symmetry group.
All other balanced colorings are said to be exotic. We analyze balanced
colorings for two closely related types of network encountered in applications:
trained Wilson networks, which occur in models of binocular rivalry, and
opinion networks, which occur in models of decision making. We give two
examples of exotic colorings which apply to both types of network, and
prove that Wilson networks with at most two learned patterns have no exotic
colorings. We discuss how exotic colorings affect the existence and 
stability of branches for bifurcations of the corresponding model ODEs. 
\end{abstract}

\section{Introduction}

We work in the `coupled cell' network formalism of~\cite{GS06, GST05, SGP03},
which should be consulted for precise definitions and proofs. Section~\ref{S:BCA} provides
a short summary. Networks
consist of {\em nodes} connected by {\em arrows} (directed edges), both of
which are partitioned into {\em types}. Each node represents a dynamical
system, and arrows indicate couplings between these systems. Identical
types determine identical dynamics or couplings. It is sometimes convenient to interpret
a node as an `internal arrow' from that node to itself. Multiple arrows and 
self-loops are permitted; indeed, they are required to simplify the theory.

Robust patterns of synchrony in networks are classified by {\em balanced colorings} of 
the nodes. A coloring assigns a color to each node, and it is balanced if nodes of the same color
have color-isomorphic input sets --- that is, the same number of arrows of each type with
tail nodes of a given color. To each coloring is associated a subspace of the state space,
called a {\em synchrony subspace},
in which the coordinates of nodes of a given color are all equal. 
This synchrony subspace is invariant for the dynamics of any ODE whose 
structure is compatible with the network architecture --- that is, for
all {\em admissible ODEs} --- if and only if the coloring is balanced. 

The existence of these universal flow-invariant subspaces has strong implications
for bifurcations of admissible ODEs, because states lying in such a subspace
have the synchrony pattern specified by the coloring, and can be found by
restricting the ODE to the subspace. Such restrictions are precisely the admissible ODEs for
the {\em quotient network}, obtained by identifying nodes with the same color and preserving
the colored input structure. Sets of synchronous nodes are often called {\em clusters} \cite{PSHMR13}, 
and the restricted ODE on the quotient network describes the dynamics of the clusters.

Synchrony for all admissible ODEs may seem a strong condition, 
but weaker forms of synchrony that persist
under small perturbations of the ODE also correspond to balanced colorings.
This has been proved for hyperbolic equilibria (\cite[Theorem 7.6]{GST05}, \cite[Theorem 6.1]{S20rigideq}), and for hyperbolic periodic states (\cite[Theorem 6.1]{GRW10}) 
with a technical assumption on the network architecture: 
see~\cite[Appendix]{S20rigideq}.
It is plausibly conjectured to be true for all 
hyperbolic periodic states~\cite[Section 10]{GS06}, and for more complex dynamic trajectories under some
kind of assumption about persistence of the underlying attractor under small perturbations.
A central issue in this case is to find a suitable definition of this type of persistence.

In networks with symmetry, every subgroup $\Sigma$ of the symmetry group $\Gamma$ defines a
balanced coloring, where colors correspond to the orbits of $\Sigma$. We call such
a coloring an {\em orbit coloring}. The synchrony
space is then the fixed-point space of $\Sigma$. However, {\em exotic} balanced colorings that
are not of this form can exist for some symmetric networks.
This was pointed out in~\cite{GNS04} for a 12-node bidirectional ring 
with $\D_{12}$ symmetry and for certain $2$-color lattice patterns. Other examples are discussed 
in~\cite{AS07, AS08}, and for planar square and hexagonal lattices 
in~\cite{S19, SGo19,SGo20}.

\begin{figure}[htb]
\centerline{
\includegraphics[height=1.8in]{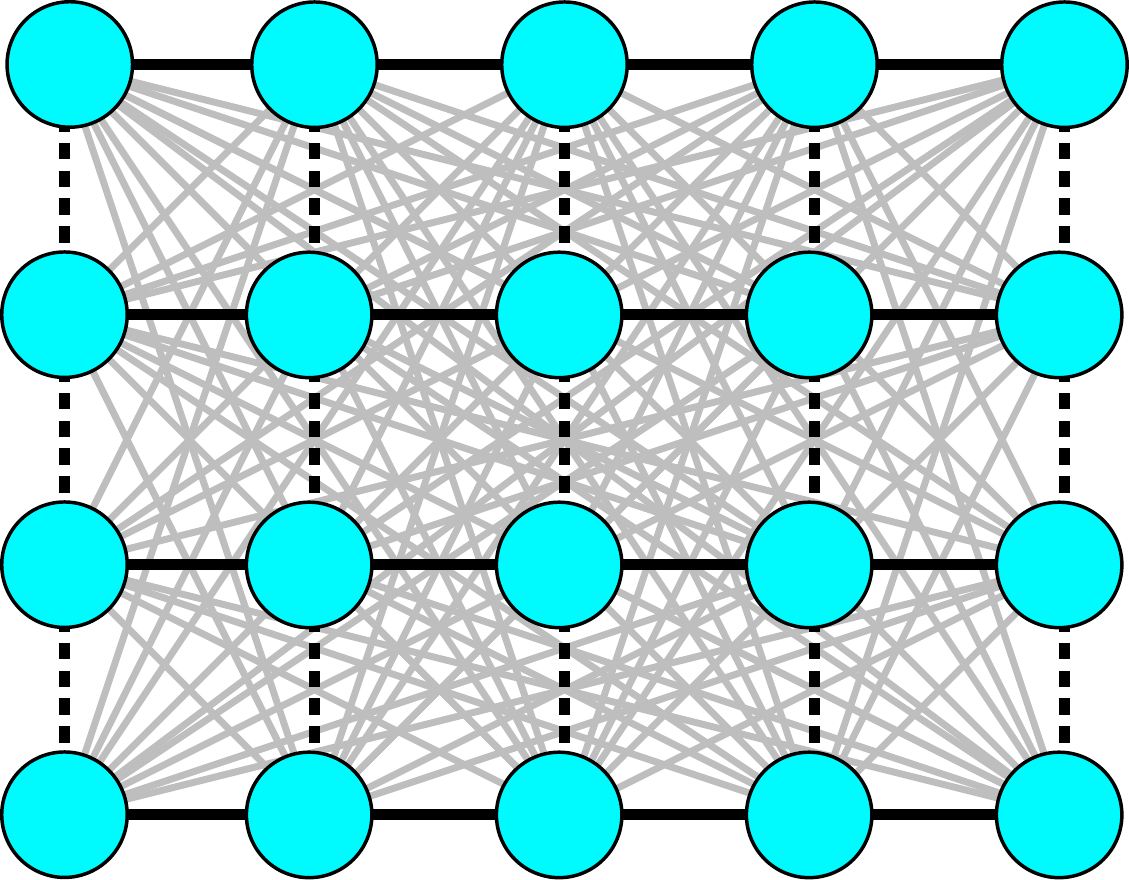}\qquad\qquad
\includegraphics[height=1.97in]{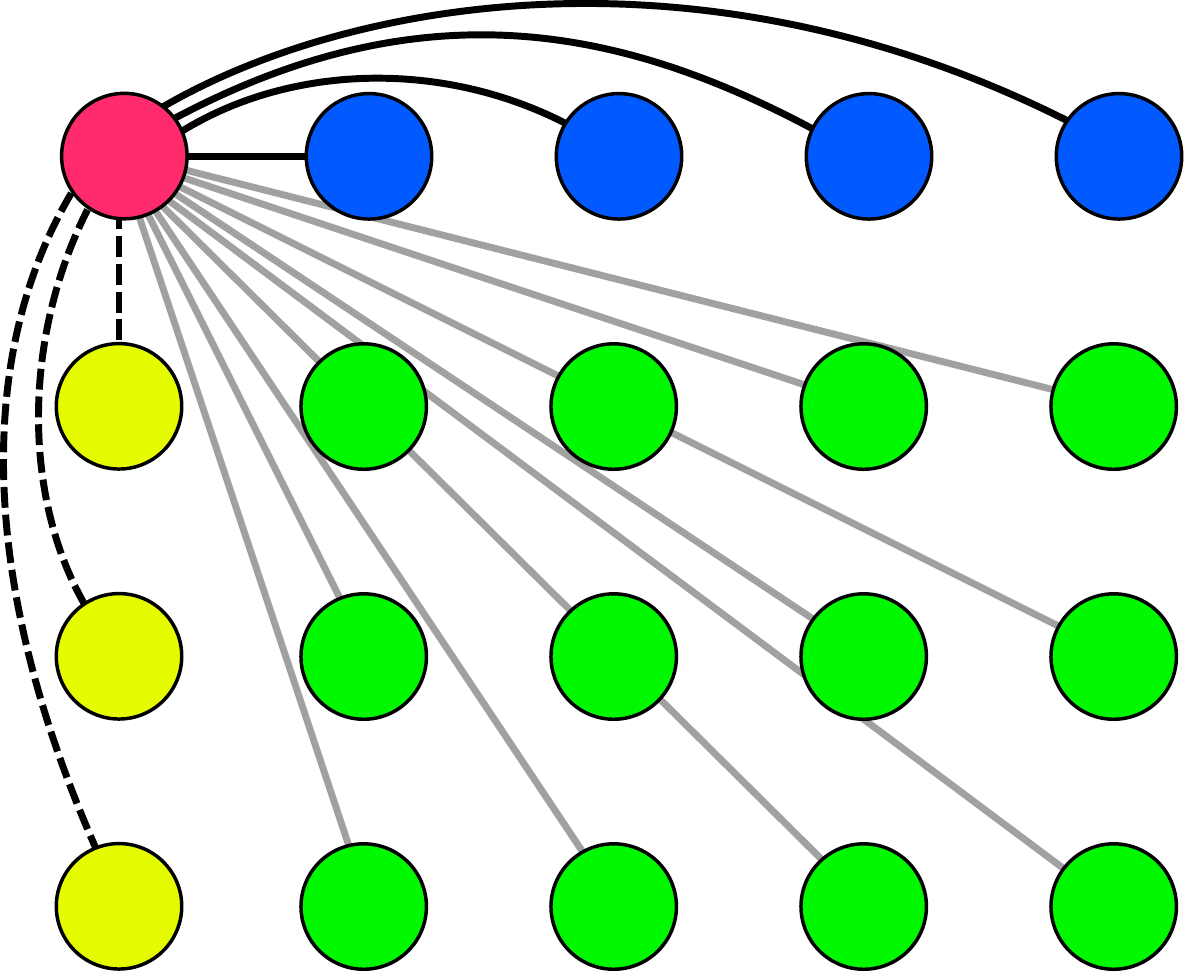}
}
\caption{{\em Left}: Group network for $\ES_4 \times \ES_5$ acting on a $4 \times 5$
array. Only a representative set of connections shown. {\em Right}:
The four arrow-types. Arrowheads not shown; only inputs to node $(1,1)$ depicted. Colors of tail nodes and forms of edges indicate arrow-type: (a) red/internal arrow, (b) blue/solid, (c) yellow/dashed, (d) green/gray.}
\label{F:5x3gpnet}
\end{figure}

Here we study networks whose nodes form an $m \times n$ array.
We assume the symmetry group is $\ES_m \times \ES_n$, where
$\ES_m$ permutes the $n$ rows and $\ES_n$ permutes the $n$ columns, so
all nodes have the same type.
Figure~\ref{F:5x3gpnet} (left) is a schematic illustration of the corresponding
{\em group network} $\GG_{mn}$ in the sense of~\cite[Section 2]{AS07}.
There is one arrow  for each ordered pair $(i,j)$ of nodes, with
head is $i$ and tail $j$. The type of the arrow is given by the orbits of such pairs
under the group action; arrows $(i,i)$ correspond to the `internal arrow' of node $i$
and are represented by the node symbol, not an arrow as such.
Figure~\ref{F:5x3gpnet} (right) is a schematic illustration of the inputs
to node $(1,1)$. There are four types of arrow:

(a)	Internal node arrow: circles.

(b)	Row arrows: solid black lines.

(c)	Column arrows: dashed black lines.

(d)	Diagonal arrows between nodes in neither the same row nor the same column: grey lines.

Admissible maps for a symmetric network are always equivariant under its symmetry group,
but the converse is false in general, \cite[Section 3.1]{AS07}. In equivariant
dynamical systems, every subspace that is flow-invariant for all equivariant ODEs is a 
fixed-point subspace~\cite[Theorem 4.1]{AS06}. The possibility
of exotic colorings can therefore 
be attributed to the difference between equivariant and admissible maps.
The synchrony subspace of an exotic balanced coloring supports a
pattern that is not an orbit coloring, hence not detected by the usual 
methods of equivariant dynamics. 

The main point of this paper is the existence of exotic colorings 
in the networks $\GG_{36}$ and $\GG_{55}$, shown below in
Figures~\ref{F:3x6wilson} and \ref{F:5x5latin} respectively. (Similar methods 
can presumably extend this result to many larger values of $m$ and $n$.)
To complement these results, we also investigate networks arising in
a model of binocular rivalry, showing that when these networks are trained
on only one or two images, no exotic colorings exist. 
We postpone discussion to Section~\ref{S:WRN},
after the models concerned have been introduced.

Finally, in Section~\ref{S:BEP} we discuss implications for bifurcations of admissible 
ODEs on the networks $\GG_{mn}$. Exotic colorings indicate the presence of
flow-invariant subspaces that are not detected by the usual method for
finding bifurcating branches of symmetric dynamical systems, namely
the Equivariant Branching Lemma. This proves the existence of
certain symmetry-breaking patterns. The analogous phenomenon of synchrony-breaking
patterns in networks can lead to exotic patterns, not predicted by their symmetries.
Taking network constraints into account, as well as
symmetry, gives a broader picture of the dynamics and bifurcations.

\section{Networks Occurring in Applications}

\begin{figure}[htb]
\centerline{%
\includegraphics[width=3in]{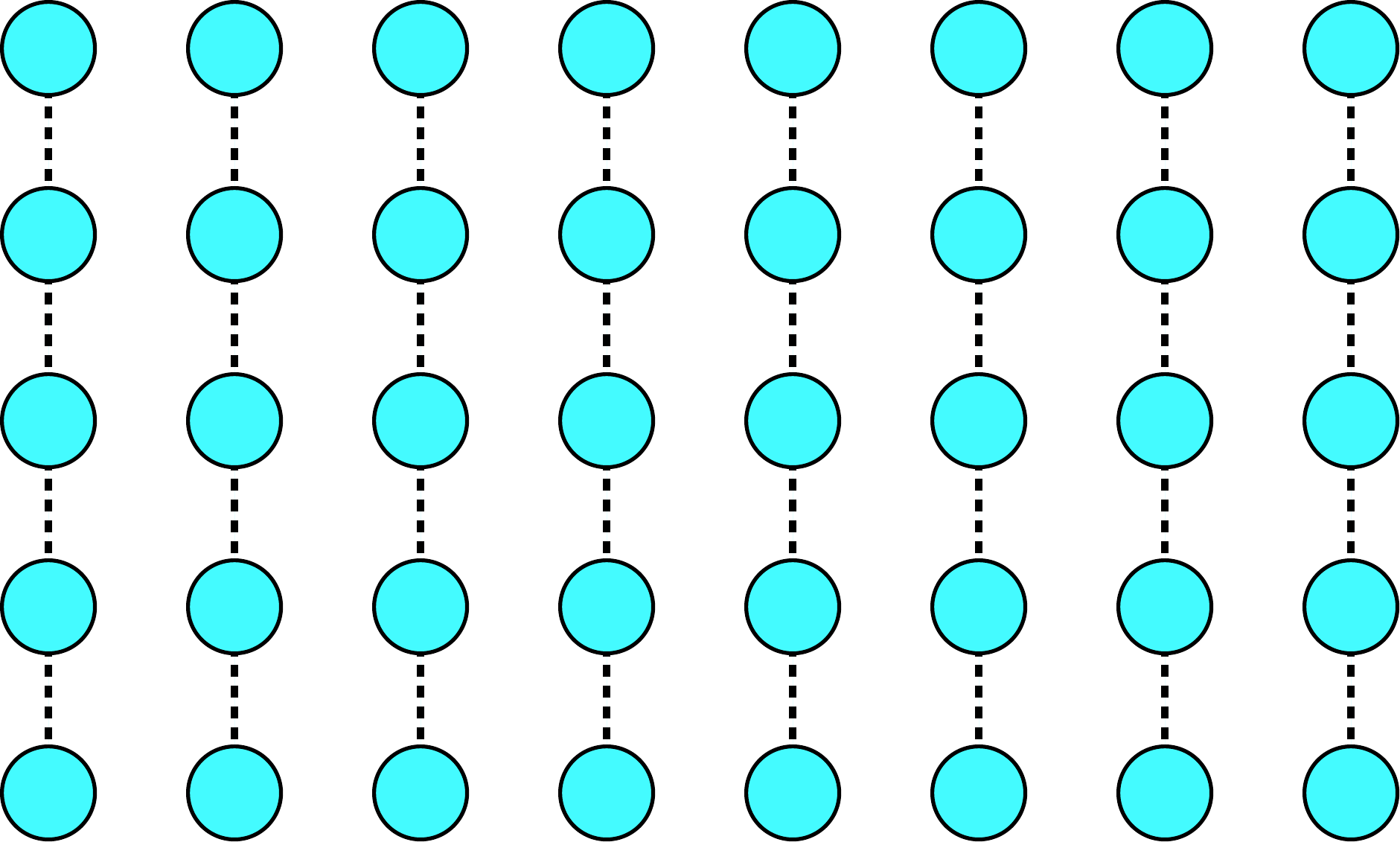}}
\caption{Untrained $5 \times 8$ Wilson network. Dotted lines indicate
sets of nodes with identical all-to-all coupling. Columns are decoupled from each other.}
\label{F:wilson_net_un}
\end{figure}

Our results apply to 
two closely related classes of symmetric networks occurring
in applications. Both are variations on networks proposed by Wilson~\cite{W03, W07, W09}
to model interocular grouping, and more generally, high-level
decision-making in the brain, and are called {\em Wilson networks}.
An {\em untrained} Wilson network is a rectangular $m\times n$ array of
identical nodes whose columns are all-to-all connected by identical inhibitory arrows. The symmetry
group of this network is not $\ES_m\times\ES_n$, but the
larger wreath product $\ES_m \wr \ES_n$ (permute the nodes in each column in any manner, and permute
the columns setwise). Figure~\ref{F:wilson_net_un} shows (schematically)
an untrained $5 \times 8$ Wilson network. For clarity, the dotted lines indicate
sets of nodes --- the columns --- with identical all-to-all coupling.

\subsection{Binocular Rivalry}

The discovery of (binocular) rivalry is credited to  Giambattista della Porta~\cite{P1593} in 1593. 
He placed two books so that each eye saw only one, and stated that he was able to
read from either book at will by moving the `visual virtue' (attention) from one eye 
to the other. Closely related but distinct phenomena are visual illusions, 
where a single image can be perceived
in more than one way; the classic example is the Necker cube of 1832, Necker~\cite{N32}, which alternates between two apparent orientations. Numerous authors have modeled
rivalry and illusions using a network with two nodes, representing the two percepts \cite{C10a,CSRR08,KM00,LC02,L88,M84,M90,NVNV07,SC11}.

Modern experiments on rivalry show a different image to each eye and ask the subject
to report what they perceive.
In simple cases, these percepts alternate between those two images.
In some studies, however,
some subjects perceive not only the original two images, 
but images that combine aspects of those images in new ways~\cite{KPYF96,SSH08,SG02,TMB06}.
The occurrence of these {\em derived images} can be modelled using trained
Wilson networks~\cite{DGMW12, DGW13, DG14}, also called rivalry networks. Experimental evidence supporting such models
has been obtained in~\cite{GZWL19}, which reports data on
the 4-location experiment proposed in~\cite{TMB06}, 
obtaining results consistent with those predicted by a trained Wilson network model. 

The model for a given pair of images begins with an untrained Wilson network.
Each column represents an
{\em attribute} of the image, such as the orientation or color of some feature. The nodes
in a given column represent {\em levels} of that attribute, such as specific orientations or colors.
The initial Wilson network is {\em trained} by adding 
excitatory arrows representing the {\em learned patterns} 
shown to the two eyes. A learned pattern is represented by a choice
of one node from each column --- the appropriate level of that attribute --- connected all-to-all by identical excitatory arrows.  

In standard model ODEs, the state of each node $i$ is defined by two variables $(x_i^A, x_i^H) \in \R^2$.
The first is an {\em activity} variable, the second a {\em fatigue} variable.
Derived patterns emerge from the dynamics of an admissible ODE for the trained Wilson network.
A balanced coloring of a Wilson network defines a {\em fusion state} in which the 
percept is ambiguous. If the fusion state becomes unstable,
symmetry-breaking states can bifurcate. 
Such bifurcations are usually analysed using equivariant (symmetric)
bifurcation theory. In the models under discussion, the relevant bifurcating branches
consist of oscillatory states, created through an equivariant version of Hopf bifurcation.
The main existence theorem here is the 
Equivariant Hopf Theorem~\cite[Chapter XVI Theorem 4.1]{GSS88},
which guarantees the existence of branches with certain spatiotemporal symmetries ---
namely, $\C$-axial subgroups, which are those having 2-dimensional fixed-point spaces
in a Liapunov-Schmidt reduction. Spatiotemporal symmetries describe patterns of
synchrony and phase relations between nodes.

The interpretation of an oscillatory state
associates the perceived image with the node or nodes whose activity variables
are largest among all nodes. Changes in the ordering of activity levels
indicate transitions between different percepts.
However, the relation of symmetry to the network structure is
subtle. In particular, exotic colorings predict behavior
in such models that is not deducible from symmetries alone. So, potentially,
exotic colorings can lead to branches with new types of spatiotemporal symmetry.

The main results of this paper for trained Wilson networks are:
\begin{itemize}
\item
Every balanced coloring is an orbit coloring for the automorphism group 
of the resulting network when the network has
been trained on 0, 1, or 2  learned patterns. 
\item
The automorphism group concerned
depends on the attribute levels of the images used to train the network.
\item There exist exotic colorings when the network has
been trained on 3 or more learned patterns.
\end{itemize}

Thus there can exist spatiotemporal patterns --- in particular, synchrony
patterns of oscillatory states and associated phase patterns --- that are not predicted by
the Equivariant Hopf Theorem. These patterns have not yet been investigated
in any detail.

\subsection{Decision Making}

The notion of decision making is very broad.
Virtually any collective system of active agents must make decisions.
The agents can, for instance, be humans, birds in a flock, fish in a school,
bees, bacteria, neurons, or autonomous robots.
As preliminary step, agents form `opinions' about a set of possible options;
that is, preferences that affect the system's behavior. These opinions collectively
determine the eventual decision at system level. Examples of
decision making systems include:

\begin{itemize}
\item
Social and economic decisions by humans based on 
different types of electoral and information-sharing systems.

\item Hunting strategies in past and present hunter-gatherer communities. 

\item Bee colonies collectively deciding on a new nest site. 

\item 
Animal groups collectively deciding when, and in which direction, to move --- for example when 
approaching two possible food sources. 

\item 
Neurons in lower brain areas integrating sensory inputs to perform perceptual and 
motor behavior decision making.

\item
Neurons in higher brain areas integrating 
sensorimotor information to make higher-level decisions.
 
\item 
Bacteria and other social microorganisms collectively deciding --- for instance, by 
quorum sensing --- how and when to undergo phenotypic differentiation in response
to environmental signals.

\item
Swarms of autonomous robots settling on a collective strategy to accomplish some task.
\end{itemize}

The same mathematical models often apply to a variety of such
real-world decision-making systems.

Here we focus on
networks very similar to Wilson networks that have been introduced into
models of decision making~\cite{FGBL20}. Again
the nodes of the network form an $m \times n$ array, but now the 
columns represent options, the rows represent agents, and the state of
a node represents the agent's opinion about that option.
Connections between nodes are of three types: distinct nodes in the same row,
distinct nodes in the same column, and all other pairs of distinct nodes.
The network has $\ES_m\times\ES_n$ symmetry where $\ES_m$ permutes
the $m$ rows and $\ES_n$ permutes the $n$ columns.

The analysis of~\cite{FGBL20} focuses on steady-state symmetry-breaking bifurcation
from a fully symmetric equilibrium, using the 
Equivariant Branching Lemma~\cite{Ci81} (see also~\cite[Lemma 1.3.1]{GS02} or \cite[Chapter XIII Theorem 3.3]{GSS88}). This result states that, subject
to certain genericity conditions, branches exist for all {\em axial} subgroups of the symmetry
group --- those with 1-dimensional fixed-point spaces. Classifying the axial subgroups
therefore leads to a list of branches breaking symmetry in various ways. In this
application the corresponding orbit colorings classify patterns of equal opinions 
about various options for various agents. These patterns govern agreements
(consensus) and disagreements (dissensus) among the agents.

The states obtained using the Equivariant Branching Lemma need not
exhaust the possible bifurcating branches, even for equivariant dynamical
systems, but they provide useful information on symmetry-breaking states
that are guaranteed to occur. In networks, the difference between
equivariant and admissible ODEs can create branches that do not occur
in general equivariant systems. In particular, exotic balanced colorings 
predict new branches that do not appear in any classification by isotropy subgroups.

The main result of this paper for opinion networks is:
\begin{itemize}
\item The examples of exotic colorings for trained Wilson networks can
also be interpreted as exotic colorings for opinion networks. 
\end{itemize}

\section{Balanced Colorings and Automorphisms}
\label{S:BCA}

We begin by recalling some basic concepts of the `coupled cell' network
formalism~\cite{GST05, SGP03}.

\begin{definition} \label{D:network_diag} \rm
A {\em network diagram} is a labelled directed graph that has four ingredients:
\begin{itemize}
\item[\rm (1)] 
A finite set of  {\em nodes} (or {\em cells}) $\CC = \{1, 2, \ldots, n\}$.
\item[\rm (2)] 
A {\em node symbol} assigned to each node. Nodes with the same symbol are said to be  {\em identical}.
\item[\rm (3)]
A finite set $\AA$ of {\em edges} or {\em arrows}. Each arrow $e$ points from its
{\em tail node} $\tl(e)$ to its {\em head node} $\hd(e)$.

\item[\rm (4)]
An {\em arrow symbol} assigned to each arrow. Arrows with the same symbol are said to be  {\em identical}.
\end{itemize}
\end{definition}

The formal theory of networks goes on to define a class of {\em admissible maps}
and associated {\em admissible ODEs} associated with the network. Essentially,
these are the maps or ODEs that respect the network structure.

The {\em automorphism group} $\aut(\GG)$ of a network $\GG$
consists of all permutations of the nodes $\CC$ that preserve the number and type
of arrows between each pair of nodes. 

The {\em input set} $I(c)$ consists of all input arrows to node $c$.
An {\em input isomorphism} $\beta:I(c) \rightarrow I(d)$ is bijection
that preserves arrow-type.

A {\em coloring} of $\GG$ is a partition of $\CC$, or equivalently
a map $K:\CC \rightarrow \KK$ where $\KK$ is a set of colors.
A coloring is {\em balanced} if, whenever nodes $c, d$ have the same color,
there is an input isomorphism $\beta:I(c) \rightarrow I(d)$ such that
$\tl(e)$ and $\tl(\beta(e))$ have the same color for all arrows $e \in \AA$.
As already remarked, a coloring is balanced if and only if it 
defines a synchrony subspace that is flow-invariant for
all admissible ODEs, and the dynamics of the resulting clusters is given by
the restricted admissible ODE on the {\em quotient network} whose nodes 
represent the clusters.

A subgroup $\Sigma \subseteq \aut(\GG)$ defines a coloring
$K^\Sigma$, where colors are determined by $\Sigma$-orbits on $\CC$.
That is, nodes have the same color if and only if they are in the same 
$\Sigma$-orbit. We call this the {\em orbit coloring} for $\Sigma$.

\begin{proposition}
Every orbit coloring is balanced.
\end{proposition}
\begin{proof}
This is~\cite[Proposition 3.3]{AS06}, stated using different terminology.
In that paper, `orbit coloring' is replaced by `fixed-point coloring'.
\end{proof}

The {\em isotropy subgroup} $\Sigma^K$ of a coloring $K$
is the subgroup of $\aut(\GG)$ that leaves every set of $K$ setwise
invariant, where $K$ is considered as a partition of $\CC$.

The next proposition is simple but useful:
\begin{proposition}
\label{P:iso_col}
A coloring is an orbit coloring if and only if it is the
same (up to choice of colors) as the orbit coloring for
$K^\Sigma$, where $\Sigma$ is its isotropy subgroup.
\end{proposition}
\begin{proof}
If $K= K^\Sigma$ then $K$ is an orbit coloring. Conversely, 
let $K$ be the orbit coloring for a subgroup $\Omega$.
Then $\Omega \subseteq \Sigma$ so $\Fix(\Sigma) \subseteq \Fix(\Omega)$.
Let $c,d$ be in the same $\Sigma$-orbit. Since $\Sigma$ fixes
the colors, $c$ and $d$ have the same color. Therefore
$\Fix(\Sigma) = \Fix(\Omega)$ and $K^\Sigma = K$.
\end{proof}

The next result is useful because some of the applications we discuss
use networks like $\GG_{mn}$ but with arrows of type (d) deleted.
\begin{lemma}
\label{L:delete}
A coloring $K$ is balanced for $\GG_{mn}$ provided it is balanced
when type (d) arrows are deleted and type (a) `arrows' are ignored.
\end{lemma}
\begin{proof}
By~\cite[Lemma 5.4]{S07} a coloring is balanced if and only if it is balanced for each type
of arrow separately --- that is, in the network where arrows of all other types are deleted.
This follows by counting arrows with tail node of a given color, and noting that
input isomorphisms preserve arrow-type. If an input isomorphism also
preserve colors of tail nodes, it must preserve colors of tail nodes of arrows of
each type separately. 

For $\GG_{mn}$, balance for type (a) internal arrows is trivial.
Moreover, balance for type (d) arrows is a consequence of balance for
both types (b) and (c), because there is a unique arrow for each pair $(i,j)$ and
it therefore has a unique type.
\end{proof}

\section{Existence of Exotic Colorings}

We now give two examples of exotic balanced colorings in networks
$\GG_{mn}$ when $(m,n) = (3,6)$ and $(5,5)$.

\subsection{The Network $\GG_{36}$}

\begin{figure}[htb]
\centerline{%
\includegraphics[width=2.7in]{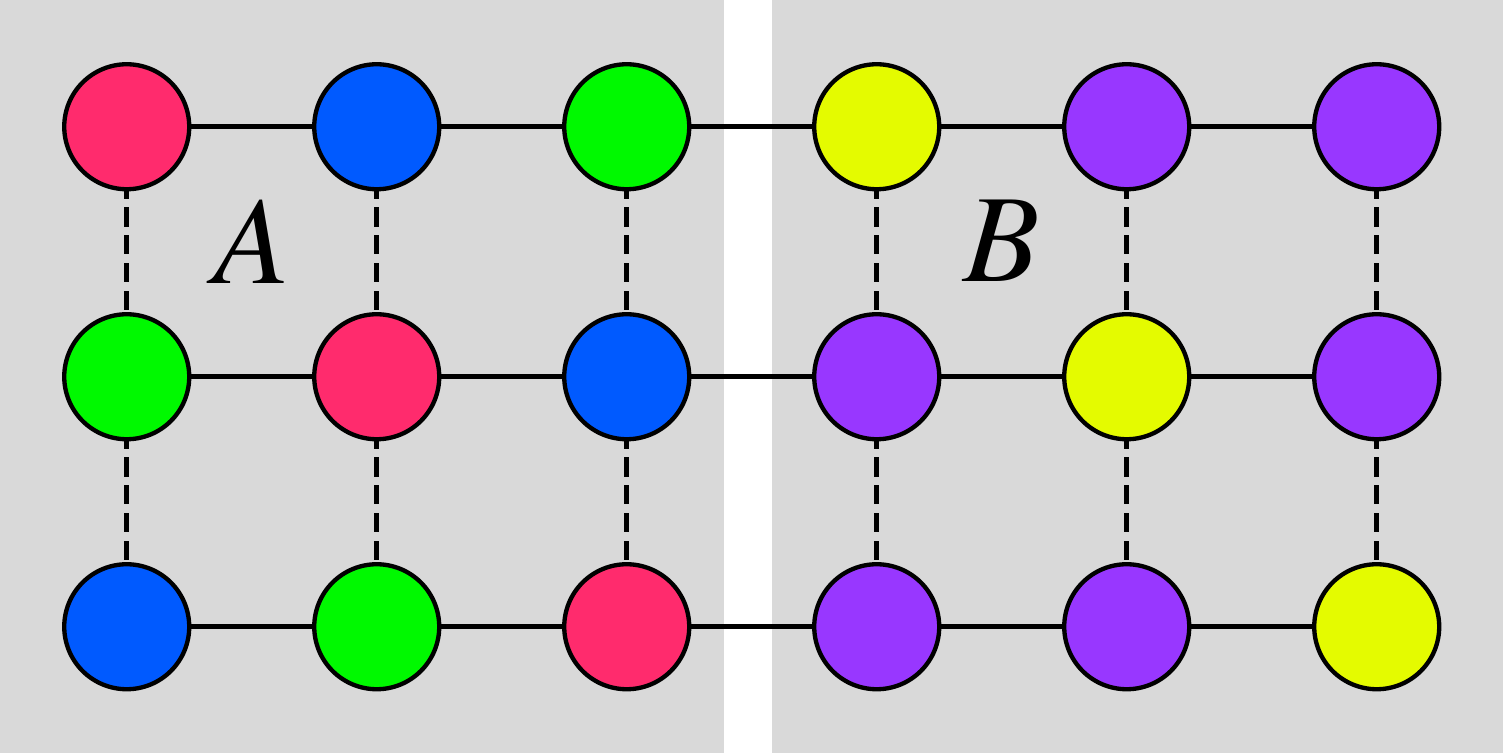}}
\caption{Exotic balanced coloring of $\GG_{36}$. Type (d) arrows omitted,
without loss of generality.}
\label{F:3x6wilson}
\end{figure}

The first example is $\GG_{36}$, Figure~\ref{F:3x6wilson}, where we
omit type (d) arrows and show types (b) and (c) schematically.
By Lemma~\ref{L:delete} the omission of type (d) arrows does not affect balance.

\begin{theorem}
The coloring in Figure~{\em \ref{F:3x6wilson}} is balanced, but is not an orbit coloring.
\end{theorem}

\begin{proof}
It is easy to prove that the pattern in Figure~\ref{F:3x6wilson} is balanced. 
Each row contains the same colors with the same multiplicities, so
it is balanced with respect to the type (b) connections.
Columns either have the same colors with the same multiplicities,
or a disjoint set of colors. So the pattern is balanced with respect to the type (c)
connections. Therefore it is balanced.

We claim it is not an orbit coloring. 
Suppose for a contradiction that it is the orbit coloring of a subgroup $\Sigma \subseteq \ES_3 \times \ES_6$. Because the colors in regions $A$ and $B$ are disjoint, we must have
\[
\Sigma \subseteq \ES_3^R \times \ES_3^A \times \ES_3^B 
\]
where $\ES_3^R$ permutes the rows, $\ES_3^A$ permutes
columns $1,2, 3$ and $\ES_3^B$ permutes columns $4,5,6$. We identify these groups
with the same group $\ES_3$ and use $R, A, B$ to indicate
which direct factor is intended.
We ask when an element $(\alpha,\beta,\gamma)  \in \ES_3^R \times \ES_3^A \times \ES_3^B$
fixes the pattern. This gives its isotropy subgroup $\Sigma$. Then we show that
the pattern is not given by the orbits of $\Sigma$.

Regions $A$ and $B$ in the figure (columns 1--3 and 4--6)
have disjoint sets of colors. The isotropy subgroup $\Sigma$
must preserve these colors, hence it must preserve regions
$A$ and $B$. Therefore $\Sigma$ is the intersection of
the isotropy subgroups  $\Sigma^A$ and  $\Sigma^B$ of the patterns in regions $A$ and $B$.
The isotropy group $\Sigma^A$ clearly consists of all elements of the form
\[
(\alpha, \alpha, \gamma) \quad \alpha \in \Z_3, \gamma \in \ES_3
\]
The isotropy group $\Sigma^B$ clearly consists of all elements of the form
\[
(\alpha, \beta, \alpha) \quad \alpha \in \ES_3, \beta \in \ES_3
\]
(This has two orbits: the set of nodes in positions $(i,i)\ (1 \leq i \leq 3)$ and
the set of nodes in positions $(i,j)\ (1 \leq i\neq j \leq 3)$. This is the 2-color
pattern in region $B$.)
Their intersection therefore consists of all elements of the form
\[
(\alpha, \alpha, \alpha) \quad \alpha \in \Z_3
\]
This group has order 3. It has six orbits, not five, so Lemma~\ref{L:delete}
implies that Figure~\ref{F:3x6wilson} is not an orbit coloring.
In fact, the orbit coloring for this group splits the six purple nodes in region $B$ 
into two sets of three, with different colors, such that the same color occurs
along each broken diagonal.
\end{proof}

\subsection{The Network $\GG_{55}$}
The example in the previous section involves two disjoint `blocks' of colors
in regions $A$ and $B$. By passing to more learned patterns --- five, to be precise --- we
can find an example with one block.
We use a coloring of $\GG_{55}$ based on a Latin square to give an example of an
exotic balanced coloring.  Again we
omit type (d) arrows and show types (b) and (c) schematically, appealing to
Lemma~\ref{L:delete}.

\begin{remark}\em
The classification of $3\times 3$ and $4 \times 4$ Latin squares~\cite{wiki}
shows that no analogous examples exists with 3 or 4 learned patterns.
\end{remark}

\begin{figure}[htb]
\centerline{%
\includegraphics[width=2in]{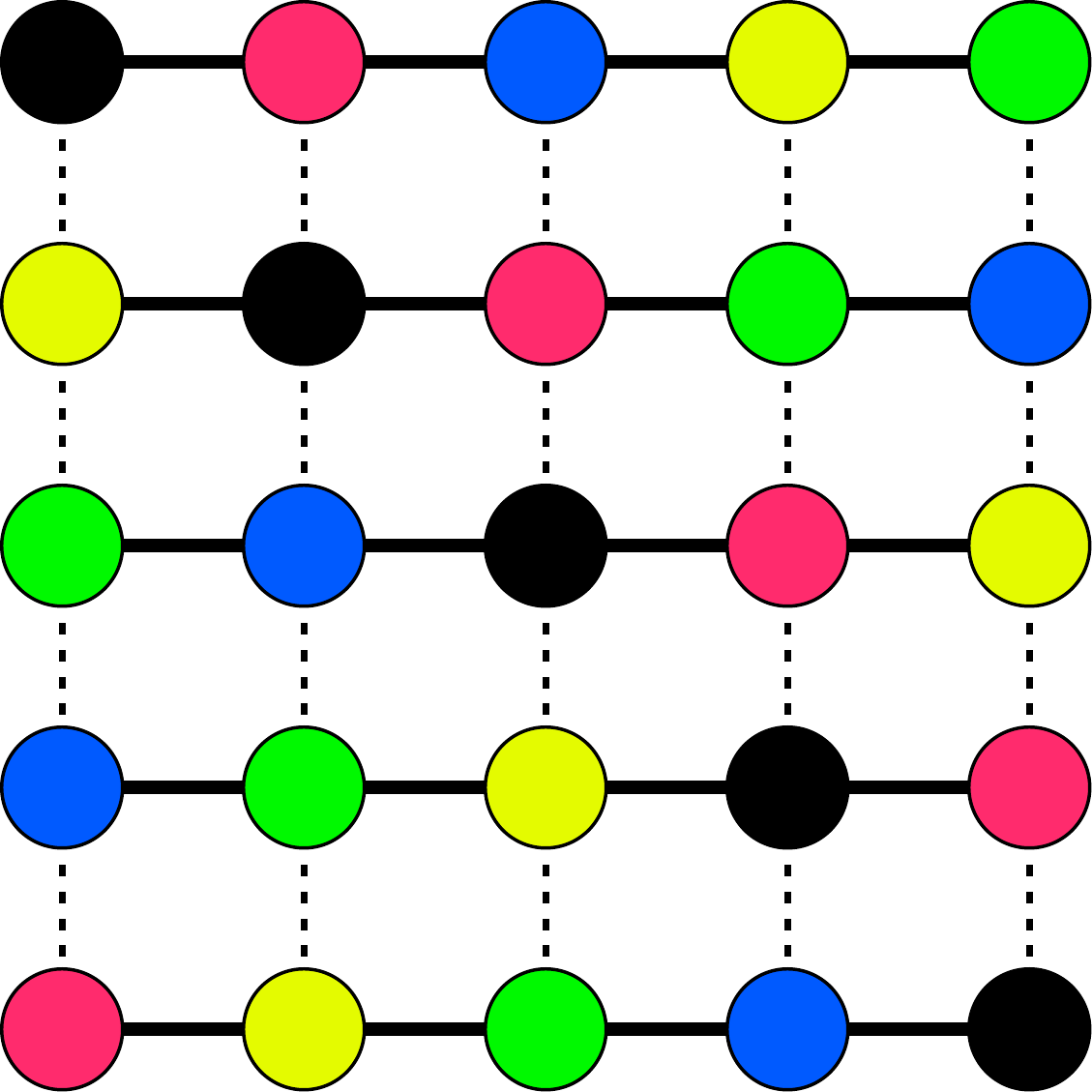}}
\caption{Exotic balanced coloring of $\GG_{55}$. Type (d) arrows omitted,
without loss of generality.}
\label{F:5x5latin}
\end{figure}

\begin{theorem}
The coloring in Figure~{\em \ref{F:5x5latin}} is balanced, but is not an orbit coloring.
\end{theorem}

\begin{proof}
This coloring is balanced because every node receives exactly one
type (b) connection and one type (c) connection from
a node of any different color.

To prove that this coloring is not an orbit coloring,
observe that the subgroup that fixes the black pattern is clearly
the diagonal subgroup
\[
\Delta = \{ (\sigma, \sigma) : \sigma \in \ES_5\} \subseteq \ES_5 \times\ES_5 
\]
For $\sigma \in \ES_5$ define the permutation matrix $P_\sigma$ by
\[
(P_\sigma)_{ij} = \left\{ \begin{array}{rcl} 
	1 & \mbox{if} & j = \sigma(i) \\
	0 & \mbox{if} & j \neq \sigma(i)
\end{array} \right.
\]

Let $M$ be a $5 \times 5$ matrix. Then $P_\sigma M$ permutes
the rows of $M$, sending row $i$ to row $\sigma^{-1}(i)$.
Similarly $MP_\sigma$ permutes
the columns of $M$, sending column $j$ to column $\sigma(j)$.
Thus an element $(\sigma, \sigma) \in \Delta$ can be considered
as acting on the space of $5 \times 5$ matrices by conjugation:
\[
M \mapsto P_\sigma MP_\sigma^{-1}
\]
This permutes both rows and columns simultaneously, via the same permutation.

Define five matrices $M_j$ corresponding to the positions of nodes of a given color:
\[
M_1 = \Matrix{1 & 0 & 0 & 0 & 0 \\ 0 & 1 & 0 & 0 & 0 \\ 0 & 0 & 1 & 0 & 0 \\
		0 & 0 & 0 & 1 & 0 \\ 0 & 0 & 0 & 0 & 1} \quad
M_2 = \Matrix{0 & 1 & 0 & 0 & 0 \\ 0 & 0 & 1 & 0 & 0 \\ 0 & 0 & 0 & 1 & 0 \\
		0 & 0 & 0 & 0 & 1 \\ 1 & 0 & 0 & 0 & 0} \quad
M_3 = \Matrix{0 & 0 & 1 & 0 & 0 \\ 0 & 0 & 0 & 0 & 1 \\ 0 & 1 & 0 & 0 & 0 \\
		1 & 0 & 0 & 0 & 0 \\ 0 & 0 & 0 & 1 & 0} 
\]
\[
M_4 = \Matrix{0 & 0 & 0 & 1 & 0 \\ 1 & 0 & 0 & 0 & 0 \\ 0 & 0 & 0 & 0 & 1 \\
		0 & 0 & 1 & 0 & 0 \\ 0 & 1 & 0 & 0 & 0} \quad
M_5 = \Matrix{0 & 0 & 0 & 0 & 1 \\ 0 & 0 & 0 & 1 & 0 \\ 1 & 0 & 0 & 0 & 0 \\
		0 & 1 & 0 & 0 & 0 \\ 0 & 0 & 1 & 0 & 0} 
\]
\noindent
We already know that the permutations $(\sigma,\tau)$ that fix the
black nodes are those in $\Delta$. So $(\sigma,\sigma) \in \Delta$ fixes
all $M_j$ (which is equivalent to fixing the coloring) if and only if 
\[
P_\sigma M_jP_\sigma^{-1} = M_j \quad 2 \leq j \leq 5
\]
or equivalently
\[
M_jP_\sigma = P_\sigma M_j \quad 2 \leq j \leq 5
\]
Solving these equations using Mathematica, we find that 
\[
P_\sigma = \Matrix{a & b & b & b & b \\ b & a & b & b & b \\ b & b & a & b & b \\
 		b & b & b & a & b \\ b & b & b & b & a }
\]
for $a, b \in \R$. This is a permutation matrix if and only if $a=1, b=0$,
so the isotropy subgroup of the coloring is the identity.
However, the fixed-point subspace of the identity contains all colorings, so
 the pattern in Figure~\ref{F:5x5latin} is not an orbit coloring.
\end{proof}

\section{Decision Making}

Networks $\GG_{mn}$ occur in models of decision making~\cite{FGBL20}.
We summarize the salient features. The model assumes that a group 
of $m$ equal agents is trying to decide among a set of $n$ options, initially
assumed to be of equal value to all agents.
This case is mathematically tractable because of its symmetry, and it is also important 
for examining more complex real-world settings. 
The equality assumptions are formalized as symmetries.
Equation (2.11) of~\cite{FGBL20} defines a specific class of admissible ODEs,
which is an equivariant
version of a more general model~\cite{Bizyaeva2020}.
The authors use the Equivariant Branching Lemma to prove the existence
of various symmetry-breaking equilibria. 

\subsection{Exotic Colorings of Opinion Networks}
\label{S:exotic_decision}

We now discuss the difference between admissible maps for the network
$\GG_{mn}$ and equivariant maps for $\ES_m\times\ES_n$. Recall
from~\cite{GST05, SGP03} that the {\em vertex group} $B(i,i)$
at a node $i$ of a network is the group of all input automorphisms
on the set $I(i)$ of all input arrows of node $i$.

The proof of \cite[Theorem 4.4]{AS06} can easily be
adapted to show that for a symmetric network,
admissible maps are the same as equivariant maps
if and only if:
\begin{itemize}
\item[\rm (a)]
The network is all-to-all connected.
\item[\rm (b)]
Every vertex symmetry extends to an element of the symmetry group. 
\end{itemize}
\noindent
When the network is $\GG_{45}$, for example, Figure~\ref{F:5x3gpnet} (right)
 shows that the vertex group of any node $c=(i,j)$ of $\GG_{mn}$ is
\[
B(c,c) = \ES_4 \times \ES_3 \times \ES_{12}
\]
Condition (a) holds in this case, but
\[
|B(c,c)| = 4!\, 3!\, 16! = 3012881743872000 \qquad |\ES_4 \times \ES_5| = 2880
\]
so condition (b) does not.
Therefore admissible maps are not the same as equivariant maps for $\GG_{45}$.

\begin{remark}\em
It can easily be shown that the linear admissible maps for $\GG_{mn}$ are 
the same as the linear equivariant maps, so the usual conditions for local bifurcation
from a fully symmetric state are the same in both cases.  However, the
detailed nonlinear structure can be different, affecting the bifurcations: see
Section~\ref{S:BEP}.
\end{remark}

The same calculation and a little extra effort proves:
\begin{theorem}
\label{T:equiv_not_admiss}
The equivariant maps for $\ES_m\times\ES_n$ are the same as the admissible maps
for $\GG_{mn}$ if and only if $(m,n) = (1,n), (2,2), (2,3), (3,2)$.
\end{theorem}
\begin{proof}
The vertex group of any node has order $v = (m-1)!(n-1)![(m-1)(n-1)]!$.
This must be less than or equal to $|\ES_m\times\ES_n| = g = m!n!$. It is easy to show that this is the case if
and only if $m = 1$, and $n$ is arbitrary, or $(m,n) = (2,2), (2,3), (2,4), (3,2), (4,2)$. We
rule out $(2,4)$ and $(4,2)$ by considering the vertex group structure, not just its order.

If $m = 1$ then $v = (n-1)!$ and we must embed $\ES_{m-1}$ in $\ES_m$, which is
possible.

If $m=2, n=2$ we must embed $\ES_1\times\ES_1\times\ES_1$ in $\ES_2\times\ES_2$, 
which is possible.

If $m=2, n=3$ or the other way round we must embed $\ES_1\times\ES_2\times\ES_2$ in $\ES_2\times\ES_2$,which is possible.

If $m=2, n=4$ or the other way round we must embed $\ES_3\times\ES_3$ in $\ES_2\times\ES_4$.
This is {\em not} possible. In $\ES_3\times\ES_3$ there exist two non-conjugate order-3
elements, but in $\ES_2\times\ES_4$ all order-3 elements are conjugate.

It remains to show that in the cases listed, all equivariants are admissible.
This is trivial for $m = 1$ since it holds for an $\ES_n$ group network.
For $(2,2)$ the vertex groups are trivial. For $(2,3)$, hence also $(3,2)$, the
vertex groups are $\Z_2 \times \Z_2$ with each factor acting independently on
two arrows of a given type (solid, gray). These elements extend to
automorphisms. The networks are all-to-all connected, so conditions (a) and (b) hold,
which completes the proof.
\end{proof}

This theorem suggests that there is scope for
exotic colorings of opinion networks to exist. We confirm this by giving
two examples of such colorings.

\begin{example}\em
In this example $m=3, n = 6$ or $m=6, n = 3$.

The $3 \times 6$ array in Figure~\ref{F:3x6wilson} can be interpreted as an opinion network with
$3$ agents and $6$ options --- or, alternatively, $6$ agents and $3$ options. The trained Wilson network includes arrow of types (a), (b), (c) above, and we have seen that it is balanced for
each type separately. The opinion network has a fourth type of arrow (d). 
By Lemma~\ref{L:delete}, the same coloring is
also balanced for type (d) arrows, so it is balanced for the opinion network.
This can also be checked directly by counting all input nodes of a given color.
\end{example}

\begin{example}\em
In the second example $m=5, n = 5$.

The $5\times 5$  network of Figure~\ref{F:5x5latin} can be interpreted as
an opinion network, and again gives an exotic coloring.
\end{example}

It seems highly likely that similar constructions provide exotic patterns
(both for trained Wilson networks and opinion networks) where
$(m,n) = (k, 2k)$ for $k \geq 4$ and $(m,n) = (k, k)$ for $k \geq 6$.
We have not investigated this question.

\section{Wilson and Rivalry Networks}
\label{S:WRN}

Recall that an untrained Wilson network $\WW$ consists of $n$ disjoint columns,
as in Figure~\ref{F:wilson_net_un}.
Each column is an all-to-all connected network (each pair of distinct nodes is connected
in each direction by an arrow) in which all arrows are equivalent
and inhibitory. Arrows in distinct columns are equivalent, and each column contains
the same number $m$ of nodes. The nodes form an $m \times n$ array, and we 
 label nodes by pairs $(i,j)$ where $1 \leq i \leq m, 1\leq j \leq n$.

\begin{lemma}
The automorphism group of $\WW$ is the wreath product
\[
\Gamma = \ES_m \wr \ES_n
\]
\end{lemma}
\proof
The  assertion is a formal statement that automorphisms can permute the nodes
in each column, independently, and also permute columns.
\qed

\subsection{Colorings of Untrained Wilson Networks}

We will prove that Wilson networks trained on 0, 1, or 2 patterns
do not have exotic colorings. To do so, we require the following definition:

\begin{definition}\em
A partition of nodes is {\em color-disjoint} if nodes of the same color
belong to the same part. That is, the colors refine the partition.
\end{definition}

We also need the following proposition, which has some interest in its own right:

\begin{proposition}
\label{P:rectangles}
Let $\WW$ be an untrained Wilson network.
Every balanced coloring of $\WW$ is the orbit coloring for
an isotropy subgroup of $\Gamma$.
\end{proposition}

\proof
All nodes of $\WW$ are input equivalent.
Consider a coloring and let its isotropy subgroup be $\Sigma$.
We claim that the coloring is given by the fixed-point set of $\Sigma$;
that is, colors are determined by the $\Sigma$-orbits.

If two columns contain a node with the same color, balance implies that every
color occurs in each column with the same multiplicity. Using permutations
of the columns, we can arrange all columns that share a color so that they are next
to each other. This defines a partition of the set of columns. Call its parts `blocks'.

\begin{figure}[htb]
\centerline{%
\includegraphics[width=3in]{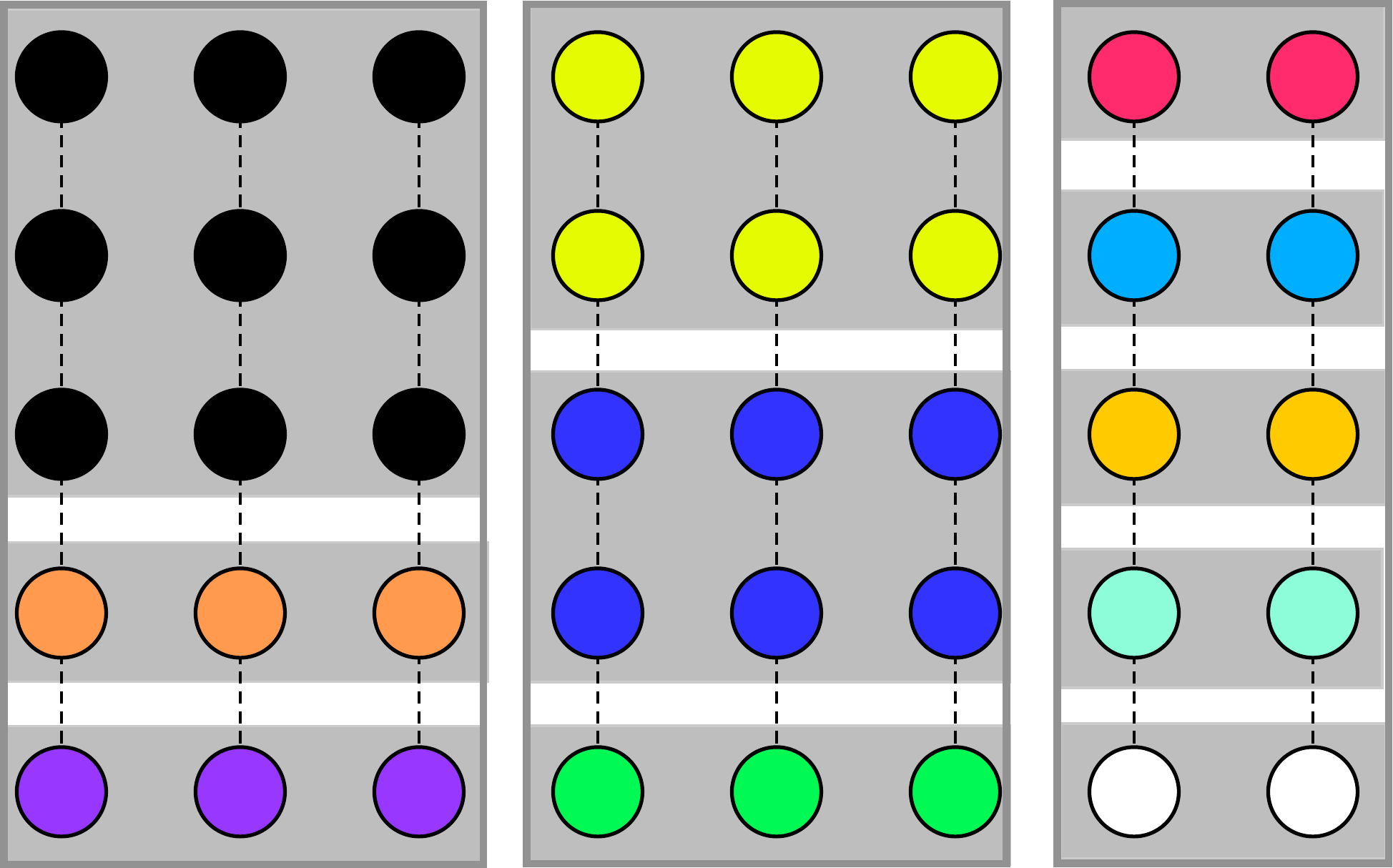}}
\caption{Partitioning an untrained Wilson network by a balanced coloring. Dotted lines indicate 
inhibitory connections within columns. Connections along each line
are all-to-all, but for clarity they are shown only for adjacent nodes.}
\label{F:wilson_net_0}
\end{figure}

Within any given block, we can permute the nodes in those columns
in any manner using $\Gamma$. Since colors occur with equal multiplicities,
we can use such permutations to arrange the colors in rectangles, 
as in Figure~\ref{F:wilson_net_0}.
Call this the {\em standard form} of the coloring. It lies in the 
$\Gamma$-orbit of the original coloring because
the permutations concerned lie in $\Gamma$. Their effect produces a coloring 
whose isotropy subgroup $\Sigma'$ is conjugate to $\Sigma$, so
without loss of generality we can consider a coloring in standard form.

We claim that the orbit coloring of $\Sigma'$
is the standard coloring. By definition the standard coloring is
fixed by $\Sigma'$. We must show that no finer coloring is fixed
by $\Sigma'$.
Since different blocks of columns 
are color-disjoint, it is enough to consider one block. 
Within that block, the rectangles are color-disjoint, 
so it is enough to consider a single rectangle.
This rectangle is a Wilson network in its own right. Since all nodes have the same color,
this is the orbit coloring for a subgroup of $\Sigma$ isomorphic to
$\ES_q \wr \ES_p$, where the rectangle has $p$ rows and $q$ columns.
\qed

\subsection{Training with Patterns}

Define
\[
M = \{ 1, \ldots, m\} \qquad
N = \{ 1, \ldots, n\} 
\]
so that nodes are labelled by $M \times N$. Here $M$ parametrizes the rows
and $N$ the columns.
A {\em pattern} is a function
\[
f:N \rightarrow M \times N
\]
such that for all $n \in N$ we have $f(m) = (g(n), n)$ for some $g(n) \in M$.
Then $g$ is a {\em section} of $M \times N$. We define $\CC^f$ to be
the set of nodes of the form $(f(n), n)$ for $n \in N$.

We represent the learning of pattern $f$ by $\WW$ as the introduction of new
excitatory arrows, leading to a network $\WW^f$. The new arrows are all
identical, and one arrow connects each ordered pair of distinct nodes in $\CC^f \times \CC^f$.
Learning a set of patterns $F=\{f_i\}$ is represented by adding one set of such
arrows for each $f_i \in F$. The resulting network, denoted by $\WW^F$, is a
{\em trained Wilson network} or {\em rivalry network}.

\begin{remark}\em
This notion of learning implies that when two (or more) learned patterns 
have some nodes in common, 
those nodes receive two (or more) new arrows, not one. 

In Wilson network models of visual illusions, learned patterns are replaced by patterns
representing pre-learned geometric consistency conditions~\cite{SG19}.
\end{remark}

We assume that all of these excitatory arrows are identical for different $f_i$.
It is easy to see that if $f$ is a single pattern then every balanced coloring
of $\WW^f$ is an orbit coloring. We prove the same is true for $2$-colorings.
(The method easily specializes to a single coloring.) First, we set
up some useful observations.

\begin{figure}[htb]
\centerline{%
\includegraphics[width=5in]{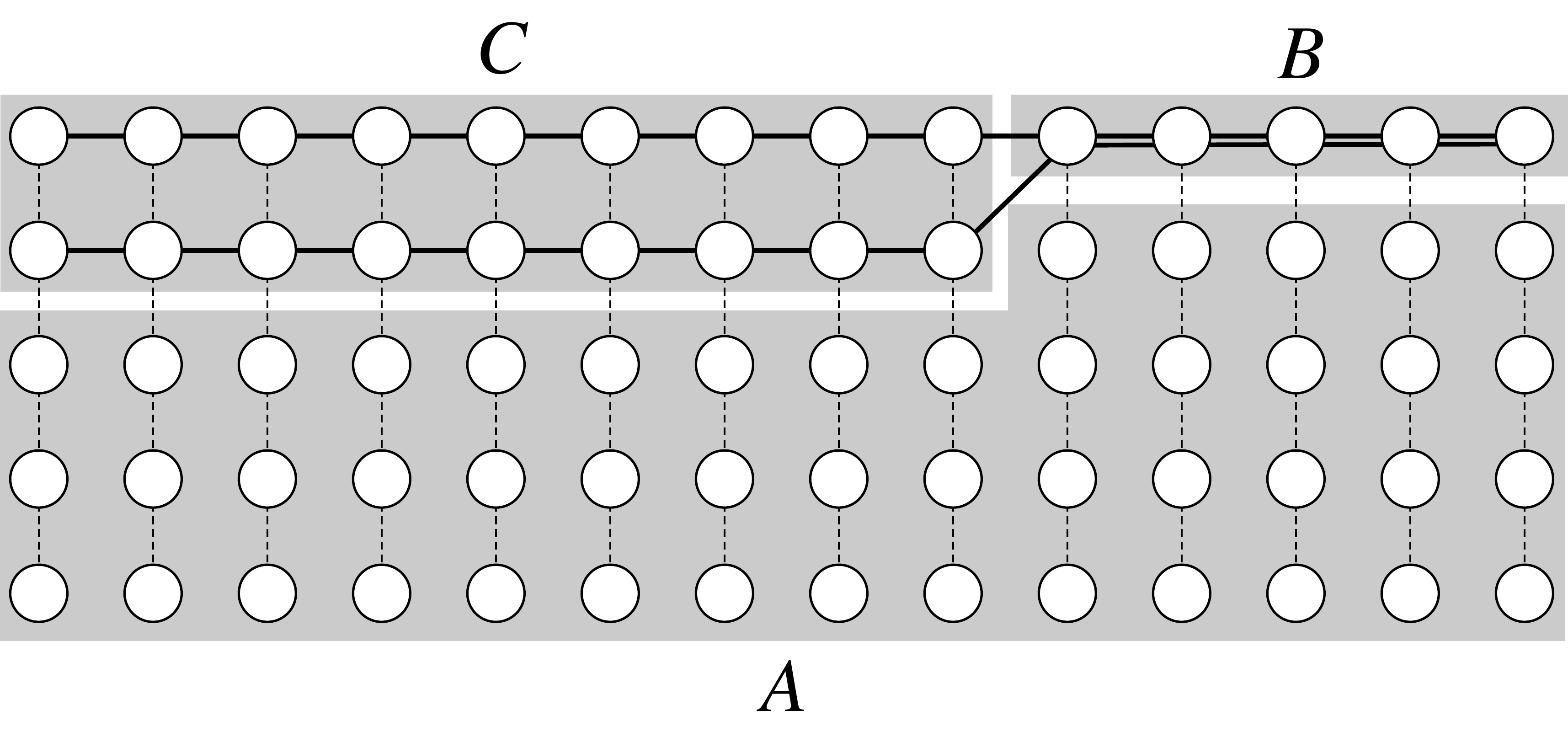}}
\caption{Partitioning a Wilson network with two learned patterns. 
Solid lines indicate learned excitatory connections,
dotted lines indicate inhibitory connections within columns. Connections along each line
are all-to-all, but for clarity these are shown only for adjacent nodes in each learned pattern.}
\label{F:wilson_net_1}
\end{figure}

Suppose that $F= \{f_1, f_2\}$ consists of two distinct patterns. Define the
{\em input-type} of a node $c$ to be the set of all nodes $d$ such that $I(d)$ is
input-isomorphic to $I(c)$. By~\cite[Section 6]{SGP03}, any
balanced coloring refines input-type. The network $\WW^F$ has precisely
three distinct input-types:
\begin{itemize}
\item[]
Type $A$: Nodes that do not appear in a learned pattern. These have no excitatory connections 
and $m-1$ inhibitory ones. 
\item[]
{Type $B$}: Nodes that appear in a single learned pattern: these have $n-1$ excitatory connections 
and $m-1$ inhibitory ones. 
\item[]
{Type $C$}: Nodes that appear in both learned patterns: these have $2n-2$ 
excitatory connections and $m-1$ inhibitory ones. 
\end{itemize}
We also use $A, B, C$ for the corresponding sets of nodes, Figure~\ref{F:wilson_net_1}.

Since a balanced coloring refines input equivalence, the partition
$\{A,B,C\}$ is color-disjoint.
These three types determine the  group $\aut(\WW^F)$.
It is a subgroup of $\aut(\WW)$ (delete excitatory connections). 
It is generated by:

\begin{itemize}
\item[\rm (1)]
Permutations of the nodes in any given column that fix all nodes in $B, C$.
These permute the nodes of type $A$ in the column.
\item[\rm (2)]
The permutation $\kappa$  that fixes every node of type $B$ and
swaps the two nodes of type $C$ in all columns that contain
such nodes.
\item[\rm (3)]
All permutations of those columns that contain nodes of type $B$.
\item[\rm (4)]
All permutations of those columns that contain nodes of type $C$.
\end{itemize}

These permutations generate the subgroup 
\begin{equation}
\label{E:Omega}
\Omega = \Z_2^\kappa \times (\ES_{m-2}\wr \ES_k) \times (\ES_{m-1}\wr \ES_{n-k}) \subseteq \ES_m \times \ES_n
\end{equation}

\begin{theorem}
\label{T:2pattorbit}
Suppose that $F$ consists of one or two distinct patterns. Then every balanced coloring
of $\WW^F$ is a orbit coloring.
\end{theorem}

\begin{proof}
We give the proof for two patterns; then we specialize to get the proof for one pattern.
So let $F = \{f_1,f_2\}$ where $f_1 \neq f_2$.

We can use permutations within each column to arrange for nodes
in learned patterns to appear as the top node, or top two nodes, of each column.
Moreover, nodes in the same $f_j$ appear in the same row 
when the patterns are distinct in that column.
Permutations between columns can position all the nodes that belong to two patterns
first (reading along rows from left to right), followed by those that belong to only one,
as in Figure~\ref{F:wilson_net_1}.

The network $\WW^F$ has precisely
three distinct input-types $A, B, C$, defined above. 
Since a balanced coloring refines input equivalence, the partition
$\{A,B,C\}$ is color-disjoint.

Suppose that the top two rows of region $C$ have a color
in common. The all-to-all excitatory connections, plus balance, then
imply that every color appearing in the top row  must also appear
in the bottom row, with the same multiplicity.

There are thus two cases:

\begin{itemize}
\item[] Case 1:
Every color appearing in the top row  must also appear
in the bottom row, with the same multiplicity.
\item[] Case 2:
The two rows of region $C$ are color-disjoint.
\end{itemize}

First, we deal with Case 1.

We refine the partition $\{A,B,C\}$ in two further stages. First we split regions $C$
into two parts $C_1, C_2$ as follows.
If a node in row 1 has an inhibitory connection to
a node with the same color in row 2, (necessarily in the same column),
place those nodes in set $C_1$. Otherwise place them in $C_2$.

We claim that each vertical connection in $C_2$ links the same
two (distinct) colors. Suppose a node in row 1 has color $\alpha$
and is connected to color $\beta$ in row 2. Balance implies
that every node of color $\alpha$ in either row must be
paired vertically with color $\beta$ in the other row. Moreover,
both rows contain the same number of nodes with any given color,
so colors $\alpha$ and $\beta$ occur with the same multiplicity
in row 1 and in row 2. (In particular, $C_2$ meets an even number of columns.)

Clearly $C_1$ and $C_2$ are color-disjoint, by balance. A color anywhere
in $C_1$ cannot connect vertically both to itself, and to something different.
(Bear in mind that regions $C$ and $A$ are color-disjoint, so the
discrepancy cannot be repaired by coloring another node in the same column.)
We now split region $A$ into $A_1, A_2, A_3$, according to whether
the top node of the column lies in $C_1, C_2,$ or $B$. Balance clearly
implies that all six regions $A_1, A_2, A_3, B, C_1, C_2$ are color-disjoint. 
See Figure~\ref{F:wilson_net_2}.

\begin{figure}[htb]
\centerline{%
\includegraphics[width=5in]{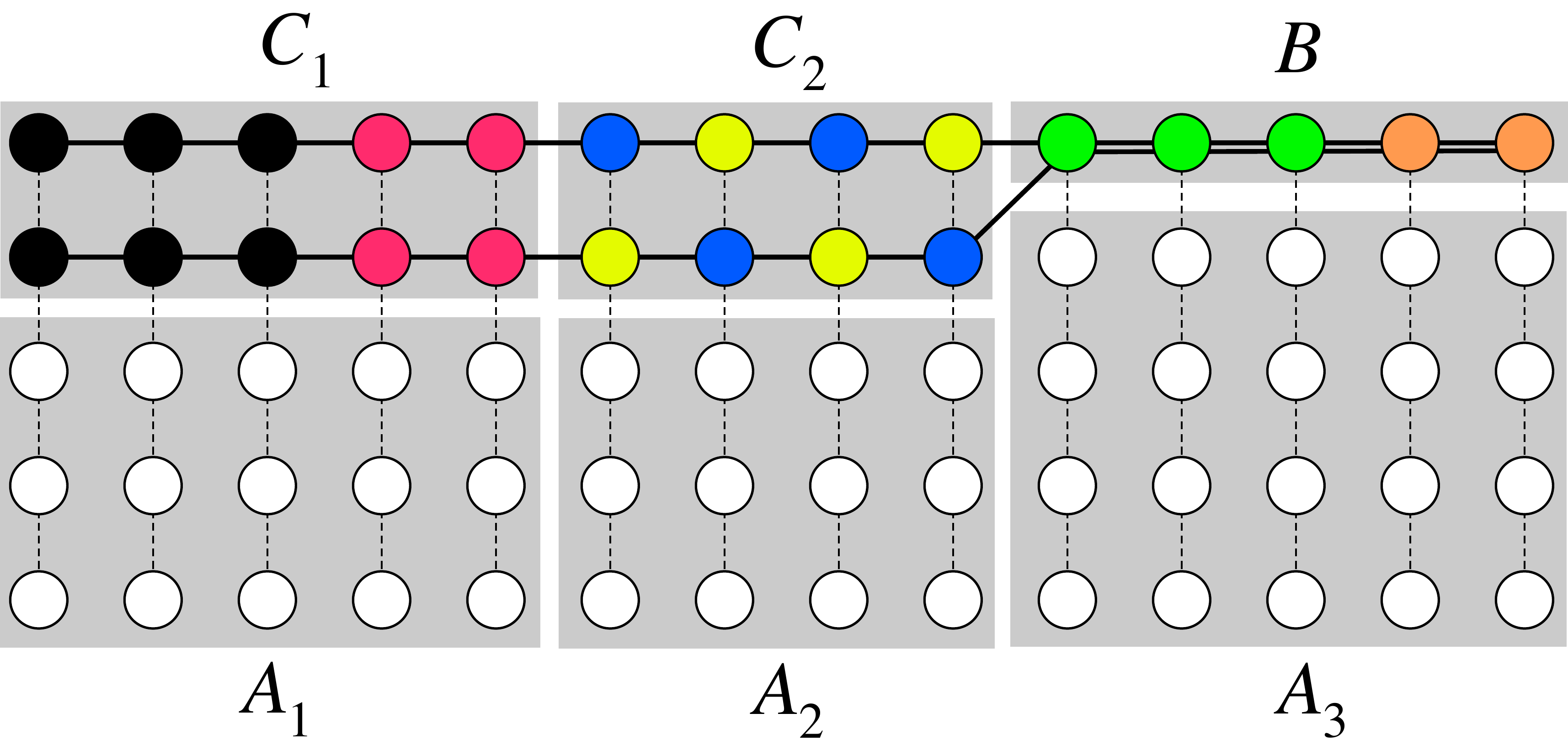}}
\caption{Refining the partition of the trained Wilson network by splitting
regions $A, C$.}
\label{F:wilson_net_2}
\end{figure}

Finally, we partition regions $C_1$ and $B$ according to the color
of the top node. Region $C_2$ is similarly partitioned using the colors
of the vertical pairs in the top two rows. (Only one such subset is
drawn but there could be many.) Without loss of generality we can
use elements of $\Omega$, 
defined in~\eqref{E:Omega}, to group nodes of a given color (or color pair)
into blocks of adjacent nodes, Figure~\ref{F:wilson_net_3}.

\begin{figure}[htb]
\centerline{%
\includegraphics[width=5in]{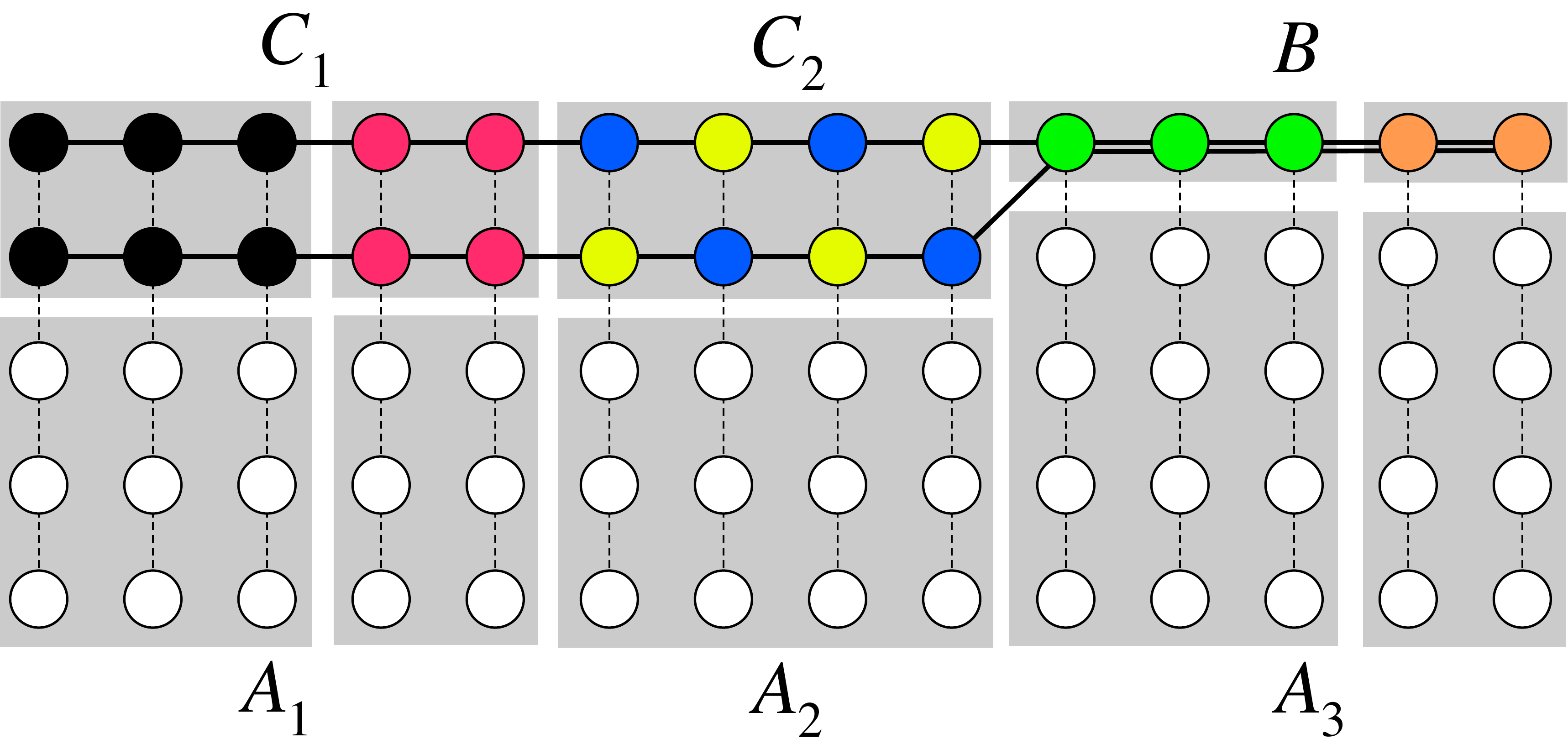}}
\caption{Refining the partition of the trained Wilson network even further
by splitting regions according to colors appearing in the learned patterns.}
\label{F:wilson_net_3}
\end{figure}

As in the proof of Proposition~\ref{P:rectangles}, we can use elements of $\Omega$
to arrange colors in the blocks of $A_j$ into rectangles, each rectangle containing
a single color.
Recall that $\kappa$ flips the two learned patterns (acting
trivially on $B$). Define a permutation $\rho$ of the columns, of order 2,
which acts trivially on $B, C_1$ and transposes pairs of adjacent columns in $C_2$
with matching colors in their top nodes. Then $\kappa\rho$ fixes all the colors
of nodes in $B \cup C$. Further permutations of the columns 
force equalities of colors along the top row.

We now construct, for any rectangle $R$ in any region, a subgroup
$\Sigma_R \subseteq \Sigma$ whose fixed-point set determines
the coloring in that rectangle. The subgroup generated by all such 
permutations, together with $\kappa\rho$ and any permutations of columns
already defined, has the same fixed-point set as $\Sigma$ because
$\Sigma$ is the isotropy group of the coloring. The construction of
$\Sigma_R$ is the same as that in the proof of Proposition~\ref{P:rectangles},
giving a group $\ES_q \wr \ES_p$ where the rectangle has $p$ rows and $q$ columns.
This completes Case 1.

Case 2 is similar, but now the two rows of $C$ are color-disjoint,
so $\kappa\rho$ is not required. 
This completes the proof for two learned patterns.

When there is only one learned pattern, the same proof applies taking $C$ to be the empty set.
Now every arrow in $B$ occurs twice, but this is effectively the same as a single
arrow of a new type.
\end{proof}

The network of Figure~\ref{F:3x6wilson} can be viewed
as a trained Wilson network with three (disjoint) learned patterns. 
Thus the analog of Theorem~\ref{T:2pattorbit} is not valid for
Wilson networks with three learned patterns. Similarly,
The network of Figure~\ref{F:5x5latin} can be viewed
as a trained Wilson network with five (disjoint) learned patterns.

\section{Bifurcations to Exotic Patterns}
\label{S:BEP}

We comment briefly on the implications of Figure~\ref{F:5x5latin} for bifurcations.
The corresponding analysis for Figure~\ref{F:3x6wilson} is less straightforward
when colors are amalgamated to apply the Equivariant Branching Lemma,
and we have not carried it out.

Associated with any balanced coloring $K$ of a network $\GG$ is a {\em quotient network} $\GG_K$
with one node for each color, and input arrows copied from the original
network according to the colors of head and tail nodes. It is shown in~\cite{GST05,SGP03}
that with canonical identifications of state spaces, the restriction of
an admissible map for $\GG$ to the synchrony subspace of $K$ is admissible for $\GG_K$,
and all admissible maps for $\GG_K$ can be obtained in this manner.

\begin{figure}[htb]
\centerline{%
\includegraphics[width=2.3in]{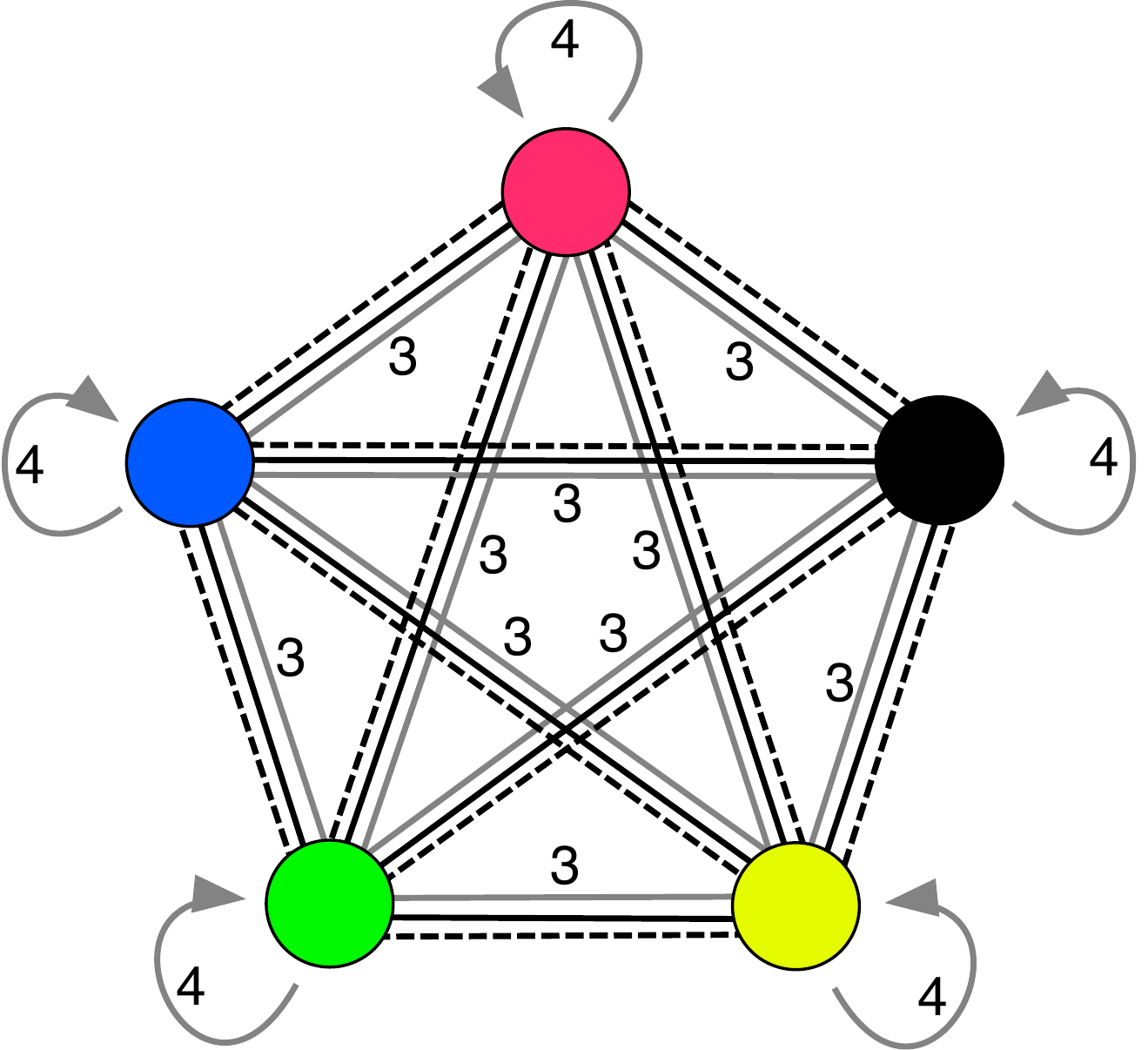}
}
\caption{The quotient network for Figure~\ref{F:5x5latin}. Lines without arrows
are bidirectional. Numbers beside gray arrows indicate multiplicities.}
\label{F:5x5latinQuot}
\end{figure}

The quotient network for the coloring in Figure~\ref{F:5x5latin} is shown in
Figure~\ref{F:5x5latinQuot}. It has five nodes and $\ES_5$ symmetry. 
It is easy to see that {\em any} coloring of this quotient network
is balanced; indeed, it is an orbit coloring for the subgroup of $\ES_5$ that preserves
the colors of nodes. 
Moreover, balanced colorings lift from quotient networks to
the original network, and remain balanced. However, the lift need not
be an orbit coloring on the full network, as we now demonstrate.

\begin{figure}[htb]
\centerline{%
\includegraphics[width=2in]{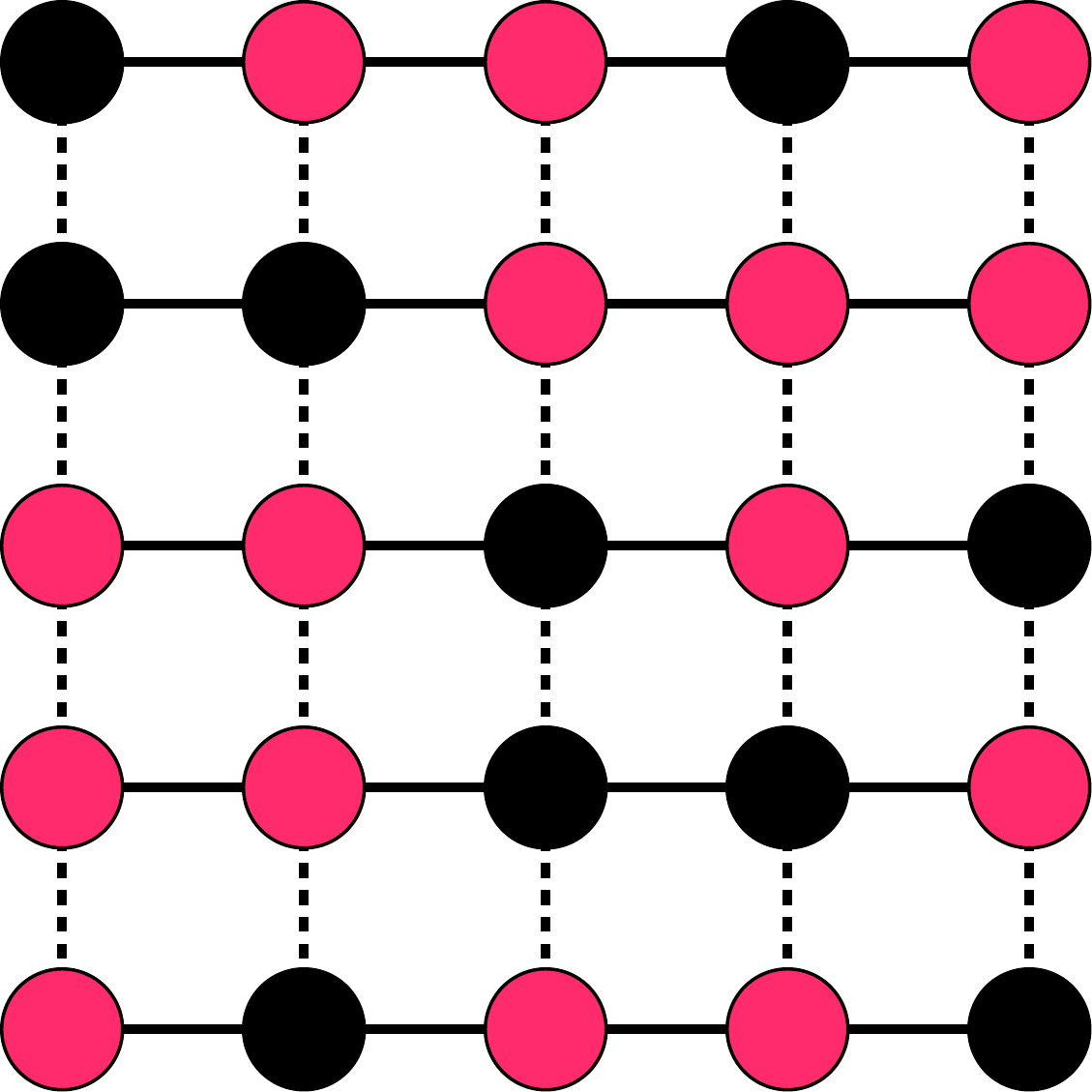} \qquad\qquad
\includegraphics[width=2in]{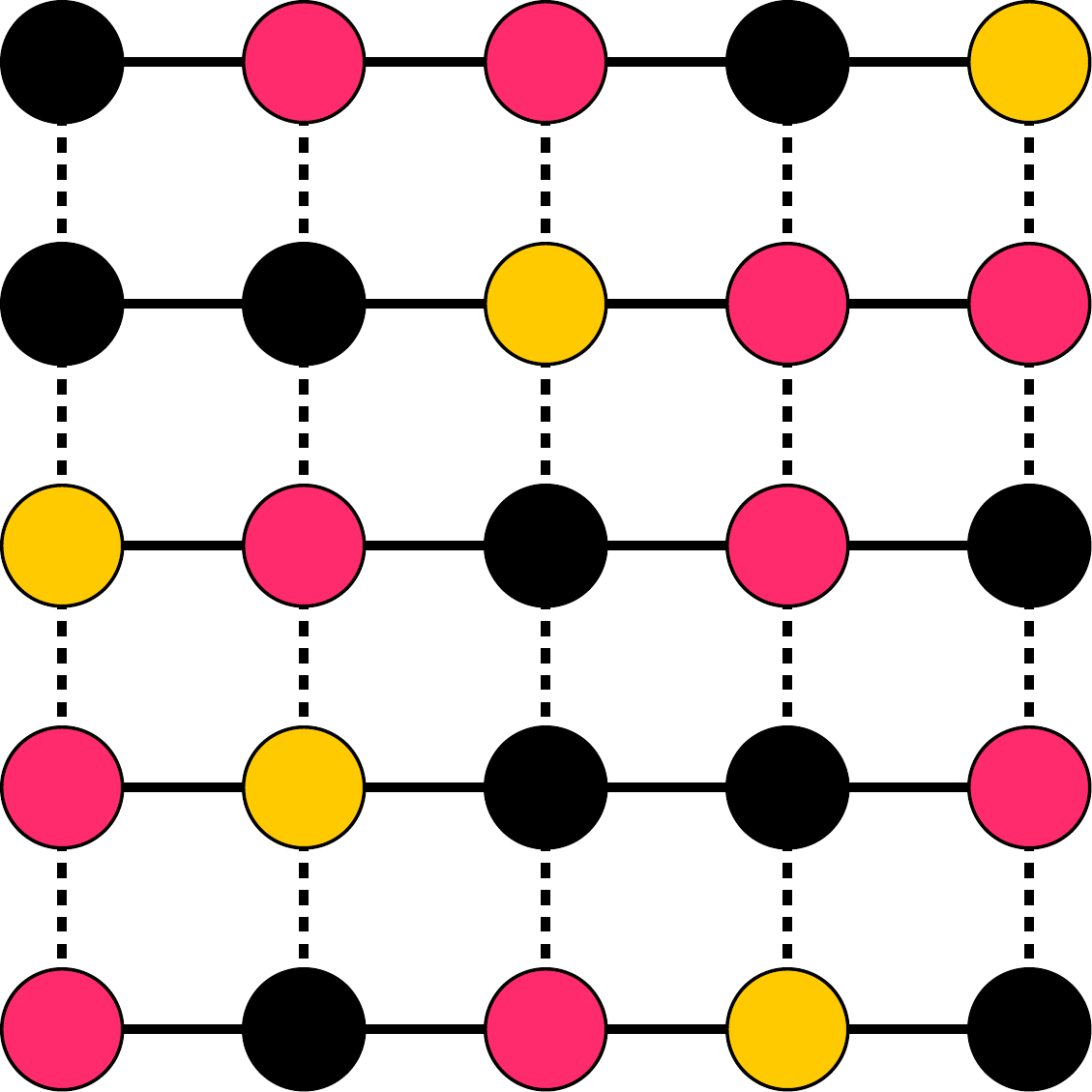}
}
\caption{{\em Left}: This balanced 2-coloring is not an orbit coloring for $\ES_5 \times \ES_5$.
{\em Right}: The fixed-point subspace of the isotropy subgroup of the 2-color pattern is a 3-color pattern.}
\label{F:5x5latin2col}
\end{figure}

We can apply the Equivariant Branching Lemma
to this quotient, and then lift the resulting pattern to $\GG_{55}$,
 to prove the generic existence of branches with (for instance) two colors, obtained
by amalgamating colors in Figure~\ref{F:5x5latin} in any manner. The question is
whether any of these amalgamated colorings is exotic for $\GG_{55}$. 
Some are not, but a plausible candidate is
Figure~\ref{F:5x5latin2col} (left), obtained by coloring black and yellow the same, 
and red, green, and blue the same. We claim that this pattern is indeed exotic
in $\GG_{55}$. The proof uses Proposition~\ref{P:iso_col}.

Calculations with Mathematica show that the
isotropy subgroup $\Sigma \subseteq \ES_5 \times \ES_5$ of
Figure~\ref{F:5x5latin2col} (left) is isomorphic
to a dihedral group $\D_5$. Since this acts on 25 nodes, and its orbits have
size at most 10, there must be at least 3 orbits, so its fixed-point subspace
must have dimension at least 3. Therefore Figure~\ref{F:5x5latin2col} is not
an orbit coloring for $\ES_5 \times \ES_5$. 

In more detail: the group $\Sigma$ is generated by 
the pairs of permutations (in cycle notation)
\[
g = ((13245),(13245)) \qquad h = ((1)(24)(35), (14)(23)(5)) 
\]
which satisfy the relations
\[
g^5 = h^2 = 1, \ hgh = g^{-1}
\]
for $\D_5$. It is a twisted $\D_5$ subgroup of $\ES_5$; that is, it consists of
elements $(\gamma, \theta(\gamma))$ where $\gamma \in \D_5 \subseteq \ES_5 \times \ONE$
and $\theta:  \D_5 \to \D_5$ is an isomorphism. Up to conjugacy we can make
$\theta$ the identity on $\Z_5$; it then maps an order-2 element to a different order-2 element.

The fixed-point subspace of $\Sigma$ comprises the patterns in 
Figure~\ref{F:5x5latin2col} (right). There are two orbits of size 10 (black, red)
and one of size 5 (orange). By Lemma~\ref{L:delete}, this $2$-coloring is not an orbit coloring.
(The orange nodes do {\em not} correspond to
a specific color in Figure~\ref{F:5x5latin} (left): the reason is that
this $5$-coloring is not an orbit coloring.)

The quotient network for the 3-coloring of Figure~\ref{F:5x5latin2col} (right),
which we do not draw, has
an obvious $\Z_2$ symmetry that swaps red and black but fixes orange.
It has an orbit coloring that merges red and black.
It also has a balanced coloring that merges red and orange, which is
not a $\Z_2$ orbit coloring, and this lifts to Figure~\ref{F:5x5latin2col} (right).

\subsection{Stability}

As well as existence of particular states, it is also important to study their stability. 
This is a more delicate issue. In particular, whenever the critical eigenspace supports
a nontrivial equivariant quadratic, all axial bifurcating branches of equilibria are {\em unstable}
near the bifurcation point; see~\cite[Chapter XIII Theorem 4.4 ]{GSS88}.
 When the symmetry group is $\ES_n$ in its standard
permutation representation, and we are considering a symmetry-breaking bifurcation, such a quadratic equivariant exists for all $n \geq 3$. This problem also arises
for $\ES_m\times\ES_n$. 
It is discussed in the evolutionary models of~\cite{CS00,GS02,SEC03},
for $\ES_n$, and in the decision models of~\cite{FGBL20} for $\ES_m\times\ES_n$. In  both cases
the remedy is to include suitable cubic equivariant terms to `compactly' the bifurcation.
This allows an unstable transcritical branch to turn around at a saddle-node (fold) point
and regain stability. This creates a jump bifurcation. 
(If $n = 2k$ there is an exception: the isotropy
subgroup $\ES_k \times\ES_k$ gives a pitchfork bifurcation.)

The same issue arises when considering the stability of the pattern in
Figure~\ref{F:5x5latin2col} (left), because the symmetry group of the quotient
network is $\ES_5$, and on the relevant synchrony subspace the pattern
is given by an axial subgroup $\ES_2 \times \ES_3$. The same remedy also
applies: include suitable cubic equivariants so that the branch turns round.
It can then regain stability within the synchrony subspace. This still leaves open
whether it is also stable transverse to the synchrony subspace, that is, to
perturbations that break the synchrony of the balanced $5$-coloring.
We have not yet investigated this question, but the problem is likely to be complicated.

Numerical simulations in~\cite{FGBL20} suggest that states arising from axial subgroups 
of $\ES_5 \times \ES_5$ can be stabilized in this manner in the whole 
of state space $\R^{25}$. However, this particular state has 
not yet been explored numerically.

\section{Conclusions}

We have studied the formation of synchrony patterns in networks $\GG_{mn}$ whose nodes form
$m \times n$ arrays, and whose connections imply $\ES_m \times \ES_n$
symmetry. Symmetry methods --- namely, the Equivariant Branching Lemma
and the Equivariant Hopf Theorem --- prove the generic existence of bifurcating
branches that break symmetry to, respectively, axial subgroups and $\C$-axial subgroups.

Network structure is stronger than equivariance, and it can create
exotic colorings, not predicted by symmetry considerations. Such colorings
exist in particular when $(m,n) = (3,6)$ and $(5,5)$.

Applying the Equivariant Branching Lemma to a quotient network of
$\GG_{55}$, we have proved the existence of an exotic 2-coloring of
$\GG_{55}$ itself. This pattern is unstable near bifurcation. It can regain
stability in the quotient network, but may still be unstable to perturbations
that break the synchrony pattern defining that quotient.

Networks $\GG_{mn}$ have been used to model decision making,
and related networks to which the same results apply occur in models
of binocular rivalry and visual illusions. The occurrence of exotic
synchrony patterns predicts states of such models that have not been
derived using equivariant methods, and whose description is not
given by isotropy subgroups.

We have also shown that in rivalry models, networks trained on one or
two images do not have exotic colorings, so all synchrony patterns 
(hence also all phase patterns) can be characterized by isotropy subgroups 
of the symmetry group.

\end{document}